\documentclass[12pt]{amsart}
\newtheorem{theorem}{Theorem}

\newtheorem{proposition}[theorem]{Proposition}
\newtheorem{corollary}[theorem]{Corollary}

\def\ds{\displaystyle}

\def\CC{\mathbb C}
\def\DD{\mathbb D}

\def\NN{\mathbb N}

\def\QQ{\mathbb Q}
\def\RR{\mathbb R}

\def\ZZ{\mathbb Z}

\def\CO{\mathcal O}

\def\CV{\mathcal V}

\DeclareSymbolFont{EulerScriptBold}{U}{eus}{b}{n}
\DeclareSymbolFontAlphabet\eusb{EulerScriptBold}

\def\Aut{\operatorname{Aut}}

\long\def\comment#1{}

\let\epsilon=\varepsilon

\let\phi=\varphi

\def\pts{\ifhmode{\leaders\hbox to 1em{\hss.\hss}\hfill}
\else{.\leaders\hbox to 1em{\hss.\hss}\hfill}\fi}

\def\Ree{\operatorname{Re}}

\def\Times_#1^#2{\overset{#2}{\underset{#1}{\bold{X}}}}
\def\too{\longrightarrow}

\def\wdtl{\widetilde}
\def\Delta{\varDelta}
\def\Gamma{\varGamma}
\def\Lambda{\varLambda}
\def\Omega{\varOmega}
\def\Phi{\varPhi}
\def\Psi{\varPsi}
\def\Xi{\varXi}

\title[Invariant metrics and distances on Neil parabolas]
{Invariant metrics and distances on generalized Neil parabolas}

\author{Nikolai Nikolov and Peter Pflug}

\address
{Institute of Mathematics and Informatics\\ Bulgarian Academy of
Sciences\\ Acad. G. Bonchev 8, 1113 Sofia,
Bulgaria}\email{nik@math.bas.bg}

\address{Carl von Ossietzky Universit\"at Oldenburg\\
Institut f\"ur Mathematik, Fakult\"at V\\ Postfach 2503\\ D-26111
Oldenburg, Germany}\email{pflug@mathematik.uni-oldenburg.de}

\subjclass[2000]{32F45}

\keywords{generalized Neil parabola, Carath\'eodory and Kobayashi
pseudodistance, Carath\'eodory-Reiffen pseudometric,
Kobayashi--Royden pseudo\-metric}

\begin{document}

\begin{thanks}{This note was written during the stay of the first
named author at the Universit\"at Oldenburg supported by a grant
from the DFG (January -- March 2006). He likes to thank both
institutions for their support.}
\end{thanks}

\begin{abstract}
We present the Carath\'eodory-Reiffen metric and the inner
Carath\'eodory distance on generalized parabolas. It turns out
that on such parabolas the Carath\'eodory distance is not inner.
\end{abstract}

\maketitle

\section{Introduction and results}

In the survey paper \cite{Jar-Pfl2} the authors had asked for an
effective formula for the Carath\'eodory distance $c_{A_{2,3}}$ on
the Neil parabola $A_{2,3}$ (in the bidisc). In a recent paper,
such a formula was presented by G.~Knese. To repeat the main
result of \cite{Kne} recall that the Neil parabola is given by
$A_{2,3}:=\{(z,w)\in\DD^2:z^2=w^3\}$, where $\DD$ denotes the open
unit disc in the complex plane. Then there is the natural
parametrization $p_{2,3}:\DD\too A_{2,3}$, $p_{2,3}(\lambda)
:=(\lambda^3,\lambda^2)$. Moreover, let $\rho$ denote the
Poincar\'e distance of the unit disc. Recall that
$\rho(\lambda,\mu)
:=\frac{1}{2}\log\frac{1+m_\DD(\lambda,\mu)}{1-m_\DD(\lambda,\mu)}$,
where
$m_\DD(\lambda,\mu):=|\frac{\lambda-\mu}{1-\lambda\overline\mu}|$,
$\lambda, \mu\in\DD$ .

Let $\lambda, \mu\in\DD$. Then Knese's result is the following
one: $$ c_{A_{2,3}}(p_{2,3}(\lambda),p_{2,3}(\mu))=
\begin{cases}
\rho(\lambda^2,\mu^2) &\text{ if } |\alpha_0|\geq 1 \\
\rho\left(\lambda^2\frac{\alpha_0-\lambda}{1-\overline\alpha_0\lambda},
\mu^2\tfrac{\alpha_0-\mu}{1-\overline\alpha_0\mu}\right) &\text{
if } |\alpha_0|<1
\end{cases},
$$ where $
\alpha_0:=\alpha_0(\lambda,\mu):=\tfrac{1}{2}(\lambda+\tfrac{1}{\overline\lambda}
+\mu+\tfrac{1}{\overline\mu}).$ In the case when $\lambda\mu=0$
the formula should be read as in the case $|\alpha_0|\geq 1$.

Observe that if $\lambda$ and $\mu$ have a non-obtuse angle, i.e.,
$\Ree(\lambda\overline\mu)\ge 0$, then $|\alpha_0(\lambda,\mu)|>1$
(compare with Corollary \ref{c=ci}).

Moreover, in \cite{Kne} the formula for the Carath\'eodory-Reiffen
pseudometric $\gamma_{A_{2,3}}$ is given as: $$
\gamma_{A_{2,3}}((a,b);X)=\begin{cases} |X_2| &\text{ if } a=b=0,
                                        |X_2|\geq 2|X_1| \\
                                        |X_1| &\text{ if } a=b=0,
                                        |X_2|<2|X_1| \\
                                        \frac{2|\lambda b|}{1-|b|^2}
                                        &\text{ if }
                                        (a,b)\neq(0,0),X=\lambda(3a,2b),
                                        \;\lambda\in\CC
                           \end{cases},
$$ where $(a,b)\in A_{2,3}$ and $X\in T_{(a,b)}A_{2,3}:=$ the
tangent space in $(a,b)$ at $A_{2,3}$.

We point out that these are the first effective formulas for the
Cara\-th\'eodory distance and the Carath\'eodory-Reiffen
pseudodistance of a non-trivial complex space.

In this paper we will discuss more general Neil parabolas, namely
the spaces $$ A_{m,n}:=\{(z,w)\in\DD^2: z^m=w^n\},\; m,n\in\NN,\
m\le n, \text{ relatively pri\-me}. $$ For short, we will call
$A_{m,n}$ the (m,n)-{\it parabola}. As in the case of the
classical Neil parabola we have the following globally bijective
holomorphic parametrization of $A_{m,n}$, namely $$
p_{m,n}:\DD\too A_{m,n},\quad
p_{m,n}(\lambda):=(\lambda^n,\lambda^m),\;\lambda\in\DD. $$
Observe that $q_{m,n}:=p_{m,n}^{-1}:A_{m,n}\too\DD$ is given
outside of the origin by $q_{m,n}(z,w)=z^kw^l$ where $k, l\in\ZZ$
are such that $kn+lm=1$; moreover, $q_{m,n}(0,0)=0$. It is clear
that $q_{m,n}$ is continuous on $A_{m,n}$ and holomorphic outside
of the origin.

We will study the Carath\'eodory and the Kobayashi distances and
also the Carath\'eodory-Reiffen and the Kobayashi-Royden
pseudometrics of $A_{m,n}$. So let us recall the objects we will
deal with in this paper: $$
m_{A_{m,n}}(\zeta,\eta):=\sup\{m_\DD(f(\zeta),f(\eta)):
f\in\CO(A_{m,n},\DD)\},\quad \zeta, \eta\in A_{m,n}, $$ where
$\CO(A_{m,n},\DD)$ denotes the family of holomorphic functions on
$A_{m,n}$, i.e., the family of those functions on $A_{m,n}$ that
are locally restriction of holomorphic functions on an open set in
$\CC^2$.

Observe that the Carath\'eodory distance $c_{A_{m,n}}$ is given by
$c_{A_{m,n}}(\zeta,\eta)\\=\tanh^{-1}m_{A_{m,n}}(\zeta,\eta)$;
moreover, $c_\DD=\rho$.

So, we have to study holomorphic function on the (m,n)-parabola.
Recall that there is the following bijection of $\CO(A_{m,n},\DD)$
and a part $\CO_{m,n}(\DD)$ of $\CO(\DD,\DD))$, where $$
\CO_{m,n}(\DD):=\{h\in\CO(\DD,\DD): h^{(s)}(0)=0, s\in S_{m,n}\}
$$ and $S_{m,n}:=\{s\in\NN:s\notin \ZZ_+ m+\ZZ_+ n\}$ (recall that
$S_{1,n}=\emptyset$ and if $m\ge2$, then $\ds\max_{s\in
S_{m,n}}s=nm-m-n$). To be precise, if $f\in\CO(A_{m,n},\DD)$ then
$f\circ p_{m,n}\in\CO_{m,n}(\DD)$, and conversely, if
$h\in\CO_{m,n}(\DD)$ then $h\circ q_{m,n}\in\CO(A_{m,n},\DD)$.

From this consideration it follows that there is the following
description of the Carathe\'odory distance on $A_{m,n}$:
\begin{align*}
&m_{A_{m,n}}(p_{m,n}(\lambda),p_{m,n}(\mu))=\max\{m_\DD(h(\lambda),
h(\mu)):h\in\CO_{m,n}(\DD)\}\\ &=\max\{m_\DD(h(\lambda), h(\mu)):
h\in\CO_{m,n}(\DD), h(0)=0\}\\
&=\max\{m_\DD(\lambda^mh(\lambda),\mu^mh(\mu)):
h\in\CO(\DD,\overline\DD),\\&\hskip 1.5cm h^{(j)}(0)=0, j+m\in
S_{m,n}\},\quad \lambda,\mu\in\DD.
\end{align*}

We like to mention that the calculation of the Carath\'eodory
distance of a generalized Neil parabola may be read as the
following interpolation problem for holomorphic functions on the
unit disc. Let $\lambda, \mu$ be as above and let $\zeta,
\eta\in\DD$. Then there exists an $h\in\CO_{m,n}(\DD)$ with
$h(\lambda)=\zeta, h(\mu)=\eta$ if and only if
$m_\DD(\zeta,\eta)\leq m_{A_{m,n}}(p_{m,n}(\lambda),
p_{m,n}(\mu))$. Note that
$m_{A_{1,n}}(p_{1,n}(\lambda),p_{1,n}(\mu))=m_\DD(\lambda,\mu)$.

From the case of domains in $\CC^n$ it is well known that the
Carath\'eo\-dory distance need not to be an inner distance (see
\cite{Jar-Pfl1}). In the case of a generalized Neil parabola it
turns out that the Carath\'eodory distance is an inner distance if
and only if $m=1$.

Recall that the associated inner distance is given by
\begin{multline*}
c_{A_{m,n}}^i(\zeta,\eta):=\inf\{L_{c_{A_{n,m}}}(\alpha): \alpha
\text{ is a } \|\cdot\|-\text{rectifiable curve in  }\\ A_{m,n}
\text{ connecting }\zeta, \eta\}, \zeta,\eta\in A_{m,n},
\end{multline*}
where $L_{c_{A_{m,n}}}$ denotes the $c_{A_{m,n}}$-length.
Obviously, $c_{A_{m,n}}\leq c_{A_{m,n}}^i$. Then we have the
following result for the inner distance.

\begin{theorem}\label{inner}
Let $\lambda,\mu\in\DD$. Then
\begin{multline*}
c_{A_{m,n}}^i(p_{m,n}(\lambda),p_{m,n}(\mu))\\ =\begin{cases}
c_\DD(\lambda^m,\mu^m)&\text{ if }\Ree (\lambda\overline\mu)\geq
\cos(\pi/m)|\lambda\mu| \\ c_\DD(\lambda^m,0)+ c_\DD(0,\mu^m)
&\text{ if }\text{otherwise}
\end{cases}.
\end{multline*}
\end{theorem}

Moreover, there is the following comparison result between the
Cara\-th\'eodory distance and its associated inner one.

\begin{corollary}\label{c=ci} Let $\lambda, \mu\in\DD$.

{\rm (a)} If
$\Ree(\lambda\overline\mu)\geq\cos(\pi/m)|\lambda\mu|$, then $$
c_{A_{m,n}}^i(p_{m,n}(\lambda),p_{m,n}(\mu))=c_{A_{m,n}}(p_{m,n}(\lambda),p_{m,n}(\mu)).
$$

{\rm (b)} If $\Ree(\lambda\overline\mu)<\cos(\pi/m)|\lambda\mu|$,
then $$
c_{A_{m,n}}^i(p_{m,n}(\lambda),p_{m,n}(\mu))=c_{A_{m,n}}((p_{m,n}(\lambda),p_{m,n}(\mu))
\text{ iff } (\lambda\overline\mu)^m<0. $$ Thus, the following
conditions are equivalent.
\begin{itemize}
\item
$c_{A_{m,n}}^i(p_{m,n}(\lambda),p_{m,n}(\mu))=c_{A_{m,n}}(p_{m,n}(\lambda),p_{m,n}(\mu))$;
\item
$c_{A_{m,n}}^i(p_{m,n}(\lambda),p_{m,n}(\mu))=c_\DD(\lambda^m,\mu^m)$;
\item $\Ree(\lambda\overline\mu)\geq\cos(\pi/m)|\lambda\mu|$ or
$(\lambda\overline\mu)^m<0$.
\end{itemize}
In particular, $c_{A_{m,n}}$ is not inner if $m>1$.
\end{corollary}

Observe that the condition
$\Ree(\lambda\overline\mu)\geq\cos(\pi/m)|\lambda\mu|$ in these
results means geometrically that $\mu$ lies inside an angular
sector around $\lambda$ of opening angle equal $\pi/m$ (compare
with Knese's result from above). Moreover, opposite to the
$A_{2,3}$-case the new area $(\lambda\overline\mu)^m<0$ (i.e., the
``rays" on which the angle between $\lambda$ and $\mu$ equals to
$\frac{(2j-1)\pi}{m},$ $j=2,\dots,m-1$) appears for $A_{m,n}$ with
$m>2.$

In order to prove Theorem \ref{inner}, we have to calculate the
Carath\'eodory-Reiffen metric $\gamma_{A_{m,n}}$ outside of the
origin.

First, let us recall its definition $$
\gamma_{A_{m,n}}((z,w);X):=\max\{|f'(z,w)X|:
f\in\CO(A_{m,n},\DD)\}, $$ where $(z,w)\in A_{m,n}$ and $X$ a
tangent vector in $(z,w)$  at $A_{m,n}$. Recall that if
$(z,w)=\zeta=p_{m,n}(\lambda)$, $\lambda\in\DD\setminus\{0\}$,
then the tangent space $T_\zeta(A_{m,n})$ at $\zeta$ is spanned by
the vector $p_{m,n}'(\lambda).$ The same holds if $m=1$ and
$\lambda=0$ whereas $T_0(A_{m,n})=\CC^2$ if $m\ge 2.$

Using the above description of $\CO(A_{m,n},\DD)$ we may
reformulate this definition in the following appropriate form
which will be used here: $$
\gamma_{A_{m,n}}(p_{m,n}(\lambda);p_{m,n}'(\lambda))=
\sup\{\frac{|h'(\lambda)|}{1-|h(\lambda)|^2}:h\in\CO_{m,n}(\DD)\}.
$$

Then we have the following result.

\begin{theorem}\label{reiffen}
Let $\lambda\in\DD.$ Then $$
\gamma_{A_{m,n}}(p_{m,n}(\lambda);p_{m,n}'(\lambda))=\frac{m|\lambda|^{m-1}}{1-|\lambda|^{2m}}.
$$
\end{theorem}

It follows from the results above (as in the case of domains in
$\CC^n$) that $\gamma_{A_{m,n}}$ is the infinitesimal form of
$c_{A_{m,n}}$ outside the origin. More precisely, if
$\lambda\in\DD\setminus\{0\},$ then
\begin{align*}
\lim_{\mu\to\lambda}
\frac{c_{A_{m,n}}(p_{m,n}(\lambda),p_{m,n}(\mu))}{|\lambda-\mu|}&=
\lim_{\mu\to\lambda}\frac{c_\DD(\lambda^m,\mu^m)}{|\lambda-\mu|}\\&=
\frac{m|\lambda|^{m-1}}{1-|\lambda|^{2m}}=
\gamma_{A_{m,n}}(p_{m,n}(\lambda);p_{m,n}'(\lambda)).
\end{align*}
Observe that the same holds if $m=1$ and $\lambda=0.$

On the other hand, note that $$
\gamma_{A_{m,n}}(0;X)=\max\{|f'(z,w)X|:
f\in\CO(A_{m,n},\DD),f(0)=0\}.$$ Then for such $f$ we have $f\circ
p_{m,n}(\zeta)=\zeta^mh(\zeta),\zeta\in\DD$, where $
h\in\CO(\DD,\overline\DD)$. Observe that $\ds\frac{\partial
f}{\partial z}(0)=\frac{h^{(n-m)}(0)}{(n-m)!}$ and
$\ds\frac{\partial f}{\partial w}(0)=h(0)$ for $m\ge 2.$ Thus, if
$X=(X_1,X_2)\in\CC^2$, then $$
\gamma_{A_{m,n}}(0;X)=\max\{|X_1\frac{h^{(n)}(0)}{n!}+X_2\frac{h^{(m)}(0)}{m!}|:
h\in\CO_{m,n}(\Bbb D),h(0)=0\}$$
$$=\max\{|X_1\frac{h^{(n-m)}(0)}{(n-m)!}+X_2h(0)|:
h\in\CO(\DD,\overline\DD),h^{(j)}(0)=0,j+m\in S_{m,n}\};$$ in
particular, $\gamma_{A_{m,n}}(0;X)=\|X\|$ if $X_1X_2=0$. Using the
first equality from above, we shall prove the following
infinitesimal result at the origin.

\begin{proposition}\label{origin} Let
$X_{\lambda,\mu}:=(\lambda^n-\mu^n,\lambda^m-\mu^m).$ Then
$$\lim_{\lambda,\mu\to 0,\lambda\neq\mu}
\frac{c_{A_{m,n}}(p_{m,n}(\lambda),p_{m,n}(\mu))}
{\gamma_{A_{m,n}}(0;X_{\lambda,\mu})}=1.$$
\end{proposition}

\begin{corollary}\label{ineq}
Let $m>1.$ Then there are points $\lambda,\mu\in\Bbb D$ such that
$$c_{A_{m,n}}(p_{m,n}(\lambda),p_{m,n}(\mu))(\lambda,\mu)>
\max\{\rho(\lambda^m,\mu^m),\rho(\lambda^{m+1},\mu^{m+1})\}.$$
\end{corollary}

It turns out that the general calculation of the
Carath\'eodory-Reiffen metric at the origin becomes much more
difficult. The next theorem may give some flavor of the nature of
this formulas.
\begin{proposition}\label{34}
Let $X=(X_1,X_2)\in\CC^2$. Then $$
\gamma_{A_{3,4}}(0;X)=\begin{cases} |X_1| &\text{ if } |X_1|\geq
2|X_2| \\
                            |X_2| &\text{ if } |X_2|\geq \sqrt{2}|X_1| \\
                            |X_1|\frac{c^3-18c+(c^2+24)^{3/2}}{108}
                            &\text{ if }
                            1<c:=2\frac{|X_2|}{|X_1|}<2\sqrt{2}
                            \end{cases} .
                            $$
\end{proposition}

It seems rather difficult to calculate an effective formula of the
Cara\-th\'e\-odory distance of $A_{m,n}$. However, we have its
value at pairs of ``opposite" points; to be more precise the
following is true.

\begin{proposition}\label{2k+1}
Let $\lambda\in\DD$, $\lambda\neq 0$. Then $$
m_{A_{2,2k+1}}(p_{2,2k+1}(\lambda),p_{2,2k+1}(-\lambda))=
\frac{2|\lambda|^{2k+1}}{1+|\lambda|^{4k+2}}. $$
\end{proposition}

Observe that now, opposite to the cases before, the number
$n=2k+1$ appears in the formula.

Finally, the discussion of the Kobayashi distance and the
Kobayashi-Royden metric on $A_{m,n}$ becomes comparably much
simpler. Let us first recall the definitions of the Lempert
function $\wdtl k_{A_{m,n}}$, the Kobayashi distance $k_{A_{m,n}}$
and the Kobayashi-Royden metric $\kappa_{A_{m,n}}$.
\begin{itemize}
\item $\wdtl k_{A_{m,n}}(\zeta,\eta):=\inf\{\rho(\lambda,\mu):
\lambda,\mu\in\DD\; \exists_{\phi\in\CO(\DD,A_{m,n})} :
\phi(\lambda)=\zeta,\;\phi(\mu)=\eta\}$, \quad $\zeta, \eta\in
A_{m,n}$;

\item $k_{A_{m,n}}:=$ the largest distance on $A_{m,n}$ below of
$\wdtl k_{A_{m,n}}$; \item $\kappa_{A_{m,n}}(\zeta;X)
:=\inf\{\alpha\in\RR_+:\exists_{\phi\in\CO(\DD,A_{m,n})}:\phi(0)=\zeta,\;\alpha\phi'(0)=X\}$,
\quad $\zeta\in A_{m,n}$, $X\in T_\zeta(A_{m,n}).$
\end{itemize}

We set $\wdtl k_{A_{m,n}}(\zeta,\eta):=\infty$ or
$\kappa_{A_{m,n}}(\zeta;X):=\infty$ if there are no respective
discs $\phi.$

Since $\CO(\DD,A_{m,n})=\{p_{m,n}\circ\psi:\psi\in\CO(\DD,\DD)\}$,
then we have the following formulas (see also
\cite{Jar-Pfl2,Kne}).

\begin{proposition}\label{Koba}
Let $\lambda,\mu\in\DD$. Then $$
k_{A_{m,n}}(p_{m,n}(\lambda),p_{m,n}(\mu))=\wdtl
k_{A_{m,n}}(p_{m,n}(\lambda),p_{m,n}(\mu))=\rho(\lambda,\mu). $$
If $\lambda\neq 0$, then
$\ds\kappa_{A_{m,n}}(p_{m,n}(\lambda);p_{m,n}'(\lambda))=
\frac{1}{1-|\lambda|^2}.$

Let $X=(X_1,X_2)\in T_{0}A_{m,n}\setminus\{0\}.$ Then $$
\kappa_{A_{m,n}}(0;X) =\begin{cases}|X_2|&\text{ if }m=1\\
\infty&\text{ if otherwise}
\end{cases}.
$$
\end{proposition}

At the end of the paper a simple reducible variety is also
discussed.

\section{Proofs and additional remarks}

We start with the proof of Theorem \ref{reiffen} which will serve
as the basic information for Theorem \ref{inner}.

\begin{proof}[Proof of Theorem \ref{reiffen}]
Recall that $$
\gamma_{A_{m,n}}(p_{m,n}(\lambda);p_{m,n}'(\lambda))=\max\{\frac{|h'(\lambda)|}{1-|h(\lambda)|^2}:
h\in\CO_{m,n}(\DD)\}. $$ Observe that if $\alpha\in\DD$ and
$\Phi_\alpha(\zeta)=\frac{\alpha-\zeta}{1-\overline\alpha\zeta},$
then $h_\alpha=\Phi_\alpha\circ h\in\CO_{m,n}(\DD)$ (use, for
example, the Fa\`a di Bruno formula) and $$
\frac{|h_\alpha'(\lambda)|}{1-|h_\alpha(\lambda)|^2}=\frac{|h'(\lambda)|}{1-|h(\lambda)|^2}.
$$ Then $$
\aligned&\gamma_{A_{m,n}}(p_{m,n}(\lambda);p'_{m,n}(\lambda))=
\max\{\frac{|h'(\lambda)|}{1-|h(\lambda)|^2}: h\in\CO_{m,n}(\DD),
h(0)=0\}\\&= \max\{\frac{|(\lambda^m\wdtl
h(\lambda))'|}{1-|\lambda^m\wdtl h(\lambda)|^2}: \wdtl
h\in\CO(\DD,\overline\DD),\wdtl h^{(j)}(0)=0, j+m\in S_{m,n}\}
\\&{\begin{aligned}=|\lambda|^{m-1}
\max\{&\frac{|mh(\lambda)+\lambda
h'(\lambda)|}{1-|\lambda^mh(\lambda)|^2}:\\
&h\in\CO(\DD,\overline\DD),h^{(j)}(0)=0,j+m\in
S_{m,n}\}=\frac{m|\lambda|^{m-1}}{1-|\lambda|^{2m}}.
\end{aligned}}\endaligned
$$ The last equality is a consequence of the fact that  the
unimodular constants are the only extremal functions for $$
\max\{\frac{|mh(\lambda)+\lambda
h'(\lambda)|}{1-|\lambda^mh(\lambda)|^2}:h\in\mathcal
O(\DD,\overline{\DD})\}. $$ To prove this fact, observe that
$(h(\lambda),h'(\lambda))$ varies on all pairs $(a,b)$ satisfying
$\ds |b|\le\frac{1-|a|^2}{1-|\lambda|^2}$. Thus, we have to show
that if $0\le c,s<1$ and $0\ds\le t\le t_s:=\frac{1-s^2}{1-c^2}$,
then $F(s,t)<F(1,0)$, where $\ds
F(s,t)=\frac{ms+ct}{1-c^{2m}s^2}$. Since $F(s,t)\le F(s,t_s)$, the
problem may be reduced to the inequality $$
\frac{m(1-c^2)s+c(1-s^2)}{1-c^{2m}s^2}<\frac{m(1-c^2)}{1-c^{2m}}
\iff\frac{c(1-c^{2m})}{m(1-c^2)}<\frac{1+c^{2m}s}{1+s}. $$ Using
the inequality $\ds\frac{1+c^{2m}}{2}<\frac{1+c^{2m}s}{1+s}$, one
has to see that $$ \frac{c(1-c^{2m})}{m(1-c^2)}<\frac{1+c^{2m}}{2}
\iff 2c\sum_{j=0}^{m-1}c^{2j}<m(1+c^{2m}). $$ Finally, by summing
up the inequalities $1-c^{2j+1}>c^{2m-2j-1}(1-c^{2j+1})$ for
$j=0,\dots,m-1$, the last inequality follows.
\end{proof}

Now, we are in the position to prove Theorem \ref{inner}.

\begin{proof}[Proof of Theorem \ref{inner}]

Set
$\Lambda_{\lambda,m}=\{\zeta\in\DD:\Ree(\lambda\overline\zeta)\geq
\cos(\pi/m)|\lambda\zeta|\}$, $\lambda\in\DD$, $m\in\NN$. Recall
again that $\Lambda_{\lambda,m}$ is an angular sector around
$\lambda$.

In a first step we shall prove that if $\lambda\in\DD$ and
$\mu\in\Lambda_{\lambda,m}$, then $$
c_{A_{m,n}}^i(p_{m,n}(\lambda),p_{m,n}(\mu))=c_\DD(\lambda^m,\mu^m).
$$ Since $$ c_{A_{m,n}}^i(p_{m,n}(\lambda),p_{m,n}(\mu))\ge
c_{A_{m,n}} (p_{m,n}(\lambda),p_{m,n}(\mu))\ge
c_\DD(\lambda^m,\mu^m),\leqno{(1)} $$ we have only to prove the
opposite inequality. After rotation, we may assume that
$\lambda\in[0,1)$. By continuity, we may also assume that
$\lambda,\mu\neq 0$ and $\arg(\mu)\in(-\pi/m,\pi/m)$. Then the
geodesic for $c^{i}_\DD(\lambda^m,\mu^m)$ does not intersect the
segment $(-1,0]$. Denote by $\alpha$ this geodesic and by
$\alpha_m$ its $m$-th root ($\root m\of 1=1$). Observe that if
$\zeta,\eta\in A_{m,n}^\ast:=A_{m,n}\setminus\{0\}$, then
\begin{multline*}
c^i_{A_{m,n}}(\zeta,\eta)=\inf\{\int_0^1
\gamma_{A_{m,n}}(\alpha(t);\alpha'(t))dt:\alpha:[0,1]\to
A_{m,n}^\ast\\ \text{ is a }C^1\text{-curve connecting
}\zeta,\eta\}
\end{multline*}
(see Theorem 4.2.7 in \cite{Kob}).

It follows by Theorem \ref{reiffen} that $$
\begin{aligned}
c^{i}_{A_{m,n}}(p_{m,n}(\lambda),&p_{m,n}(\mu))\le\int_0^1
\gamma_{A_{m,n}}(p_{m,n}\circ\alpha_m(t);(p_{m,n}\circ\alpha_m)'(t))dt
\\ &=\int_0^1 \frac{m|(\alpha_m(t))|^{m-1}\alpha'_m(t)|}
{1-|\alpha_m(t)|^{2m}}dt
=\int_0^1\frac{|\alpha'(t)|}{1-|\alpha(t)|^2}dt\\ &=c^{i}_\DD
(\lambda^m,\mu^m)=c_\DD(\lambda^m,\mu^m).
\end{aligned}
$$

It remains to prove that if $\mu\not\in\Lambda_{\alpha,m}$, then
$$ c^{i}_{A_{m,n}}(p_{m,n}(\lambda),p_{m,n}(\mu))=
c^{i}_{A_{m,n}}(p_{m,n}(\lambda),0)+
c^{i}_{A_{m,n}}(0,p_{m,n}(\mu)). $$ By the triangle inequality, we
only have to prove  that $$
c^{i}_{A_{m,n}}(p_{m,n}(\lambda),p_{m,n}(\mu))\ge
c^{i}_{A_{m,n}}(p_{m,n}(\lambda),0)+
c^{i}_{A_{m,n}}(0,p_{m,n}(\mu)).\leqno{(2)} $$ Take an arbitrary
$C^1$-curve $\alpha:[0,1]\to A_{m,n}^\ast$ with $\alpha(0)=
p_{m,n}(\lambda)$ and $\alpha(1)=p_{m,n}(\mu)$. Let $t_0\in(0,1)$
be the smallest numbers such that
$\alpha(t_0)\in\partial\Lambda_{\alpha,m}.$ If $
\alpha(t_0)=p(\lambda_0),$ then $$\begin{aligned}
&\int_0^1\gamma_{A_{m,n}}(\alpha(t);\alpha'(t))dt\\
&=\int_0^{t_0}\gamma_{A_{m,n}}(\alpha(t);\alpha'(t))dt+
\int_{t_0}^1\gamma_{A_{m,n}}(\alpha(t);\alpha'(t))dt\\ &\ge
c^{i}_{A_{m,n}}(p_{m,n}(\lambda),p_{m,n}(\lambda_0))
+c^{i}_{A_{m,n}}(p_{m,n}(\lambda_0),p_{m,n}(\mu))\\ &\ge
c_{A_{m,n}}(p_{m,n}(\lambda),p_{m,n}(\lambda_0))
+c_{A_{m,n}}(p_{m,n}(\lambda_0),p_{m,n}(\mu))\\ &\ge
c_\DD(\lambda^m,\lambda_0^m)+c_\DD(\lambda_0^m,\mu^m)\\
&=c_\DD(\lambda^m,0)+c_\DD(0,\lambda_0^m)+c_\DD(\lambda_0^m,\mu^m)\quad
(\text{since }\lambda_0^m\in(-1,0)\\ &\ge c_\DD(\lambda^m,0)+
c_\DD(0,\mu^m).
\end{aligned}$$
Now, (2) follows by taking the infimum over all curves under
consideration.
\end{proof}

Next, the proof of Corollary \ref{c=ci} will be given.

\begin{proof}[Proof of Corollary \ref{c=ci}]
(a) Follows by Theorem \ref{inner} and the inequality (1).

(b) The inequalities $$ c_{A_{m,n}}(p_{m,n}(\lambda),p_{m,n}(\mu))
\le\max\{c_\DD(\lambda^mf(\lambda),\mu^mf(\mu)):f\in\CO(\DD,
\overline\DD)\} $$ $$ \le\max\{c_\DD(\lambda^mf(\lambda),0)+c_\DD
(0,\mu^m+f(\mu)):f\in\mathcal\CO(\DD,\overline\DD)\} $$ $$
\le
c_\DD(\lambda^m,0)+c_\DD(0,\mu^m)=c^{i}_{A_{m,n}}(p_{m,n}(\lambda),p_{m,n}(\mu))
$$ show that $$
c_{A_{m,n}}(p_{m,n}(\lambda),p_{m,n}(\mu))=c^{i}_{A_{m,n}}(p_{m,n}(\lambda),p_{m,n}(\mu))
$$ if and only if $\lambda^mf(\lambda)$ and $\mu^mf(\mu)$ lie on
opposite rays and $|f(\lambda)|=|f(\mu)|=1$ for some
$f\in\CO(\DD,\overline\DD)$, i.e., $f$ is a unimodular constant,
and $(\lambda\overline{\mu})^m<0$.

The remaining part of Corollary \ref{c=ci} follows by the fact
that $c_{\DD}(z,0)+c_{\DD}(0,w)=c_{\DD}(z,w)$ if and only if
$z\overline w\le 0$.
\end{proof}

\noindent{\bf Remarks.} (a) For $m\in\NN$, consider the following
distance on $\DD:$ $$ \rho^{(m)}(\lambda,\mu):=
\max\{\rho_\DD(\lambda^mh(\lambda),\mu^mh(\mu))
:h\in\CO(\DD,\overline\DD)\}. $$ Note that
\begin{align*}\lim_{\varepsilon\to
0,\varepsilon\neq
0}\frac{\rho^{(m)}(\lambda,\lambda+\varepsilon)}{|\varepsilon|}&=
|\lambda|^{m-1}\max\{\frac{|mh(\lambda)+\lambda
h'(\lambda)|}{1-|\lambda^mh(\lambda)|^2}:
h\in\CO(\DD,\overline\DD)\}\\&=\gamma_{A_{m,n}}(p_{m,n}(\lambda);p_{m,n}'(\lambda))
\end{align*}
 by the proof of Theorem \ref{reiffen}. So it follows
that the associated inner distance of $\rho^{(m)}$ equals
$c^i_{A_{m,n}}(p_{m,n}(\cdot),p_{m,n}(\cdot))$. Then
\begin{multline*}
c^i_{A_{m,n}}(p_{m,n}(\lambda),p_{m,n}(\mu))\ge
\rho^{(m)}(\lambda,\mu)\\ \ge
c_{A_{m,n}}(p_{m,n}(\lambda),p_{m,n}(\mu))\ge
\rho(\lambda^m,\mu^m).
\end{multline*}

Moreover, the proof of Corollary \ref{c=ci} shows that the
following conditions are equivalent:
\begin{itemize}
\item $c_{A_{m,n}}^i(p_{m,n}(\lambda),p_{m,n}(\mu))=
\rho^{(m)}(\lambda,\mu)$;
\item$c_{A_{m,n}}^i(p_{m,n}(\lambda),p_{m,n}(\mu))=c_{A_{m,n}}(p_{m,n}(\lambda),p_{m,n}(\mu))$;
\item
$c_{A_{m,n}}^i(p_{m,n}(\lambda),p_{m,n}(\mu))=\rho(\lambda^m,\mu^m)$;
\item $\Ree(\lambda\overline\mu)\geq\cos(\pi/m)|\lambda\mu|$ or
$(\lambda\overline\mu)^m<0$.
\end{itemize}

As an application of these observations we obtain a simple proof
(without calculations) of Lemma 14 in \cite{Pfl-Zwo}:
\smallskip

{\it If $a,b\in[0,1),$ $s\in(0,1]$ and $\theta\in[-\pi,\pi],$ then
$\rho(a,be^{i\theta})\le\rho(a^s,b^se^{is\theta}).$}\smallskip

In fact, we may assume that $s\in\QQ.$ If $s=\frac{p}{q}$ ($1\le
p\le q$), $\lambda=\root q\of a,$ $\mu=\root q\of b
e^{i\theta/q},$ then we have to prove that
$\rho(\lambda^q,\mu^q)\le\rho(\lambda^p,\mu^p).$ But the angle
between $\lambda$ and $\mu$ does not exceed
$\frac{\pi}{q}\le\frac{\pi}{p}$ and hence
$$\rho(\lambda^p,\mu^p)=\rho^{(p)}(\lambda,\mu)\ge\rho(\lambda^q,\mu^q)$$
(the last inequality holds for any $\lambda,\mu\in\DD$ and $q\ge
p$).
\medskip

(b) Recall that
\begin{align*}c_{A_{m,n}}(p_{m,n}(\lambda),p_{m,n}(\mu))
&=\max\{\rho_\DD(\lambda^mh(\lambda),\mu^mh(\mu)):\\&
h\in\CO(\DD,\overline\DD),h^{(j)}(0)=0, j+m\in S_{m,n}\}.
\end{align*}
If $m=1$ or $(m,n)=(2,3),$ then $\rho^{(m)}(\lambda,\mu)=
c_{A_{m,n}}(p_{m,n}(\lambda),p_{m,n}(\mu)),$ since
$S_{1,n}=\emptyset$ and $S_{2,3}=\{1\}.$

On the other hand, if $m\neq 1$ and $m\neq n-1$, then the
following conditions are equivalent:
\begin{itemize}

\item $\rho^{(m)}(\lambda,\mu)=\rho(\lambda^m,\mu^m);$ \item
$\rho^{(m)}(\lambda,\mu)=c_{A_{m,n}}(p_{m,n}(\lambda),p_{m,n}(\mu))$.
\end{itemize}

It is clear that the first condition implies the second one. For
the converse, observe that as $h$ varies over $\CO(\DD,\DD)$, the
pair $(h(\lambda),h(\mu))$ varies over all $(z,w)\in\DD^2$ with
$m_\DD(z,w)\le m_\DD(\lambda,\mu)$. Thus,
\begin{align*}
\rho^{(m)}(\lambda,\mu)=\max
&\{\rho_\DD(\lambda^mz,\mu^mw):z,w\in\DD \text{ with }\\
&m_\DD(z,w)\le m_\DD(\lambda,\mu)\text{ or }z=w\in\partial D.\}
\end{align*}
It follows by the maximum principle for the continuous
plurisubharmonic function $m_\DD(\lambda^m\;\cdot,\mu^m w)$ that
if $\rho^{(m)}(\lambda,\mu)= \rho_\DD(\lambda^mz,\mu^mw)$, then
either $z=w\in\partial D$, or $m_\DD(z,w)=m_\DD(\lambda,\mu)$.
Assuming that $\rho^{(m)}(\lambda,\mu)\neq\rho(\lambda^m,\mu^m)$
excludes the first possibility. Then any extremal function $h$ for
$\rho^{(m)}(\lambda,\mu)$ satisfies
$m_\DD(h(\lambda),h(\mu))=m_\DD(\lambda,\mu)$, i.e.,
$h\in\Aut(\DD)$. Since any such function should be also extremal
for $c_{A_{m,n}}(p_{m,n}(\lambda),p_{m,n}(\mu))$, it follows that
either $h^{(j)}\neq 0$ for any $j\in\NN$, or  $h$ is a rotation.
In particular, $m+1\not\in S_{m,n},$ i.e., $m=1$ or $m=n-1$, a
contradiction.

Let $m\ge 3.$ Then $m+2\not\in S_{m,m+1}$ and hence $h$ must be a
rotation. Thus, the following conditions are equivalent:

\begin{itemize}
\item
$\rho^{(m)}(\lambda,\mu)=\max\{\rho(\lambda^m,\mu^m),\rho(\lambda^{m+1},\mu^{m+1})\};$
\item
$\rho^{(m)}(\lambda,\mu)=c_{A_{m,m+1}}(p_{m,n}(\lambda),p_{m,n}(\mu))$.
\end{itemize}

(c) Concerning the first condition from above, we point out that
if $m>1,$ then by Corollary \ref{ineq} there are points
$\lambda,\mu\in\Bbb D$ such that
\begin{align*}\rho^{(m)}(\lambda,\mu)&\ge
c_{A_{m,n}}(p_{m,n}(\lambda),p_{m,n}(\mu))(\lambda,\mu)\\
&>\max\{\rho(\lambda^m,\mu^m),\rho(\lambda^{m+1},\mu^{m+1})\}.
\end{align*}
On the other hand,
$\rho^{(2m)}(\lambda,-\lambda)=\rho(\lambda^{2m+1},-\lambda^{2m+1}),
$ since $$
\begin{aligned}
m_\DD(\lambda^{2m}&\Phi_\alpha(\lambda),
\lambda^{2m}\Phi_{\alpha}(-\lambda))\\
&=\frac{2(1-|\alpha|^2)|\lambda|^{2m+1}}
{|1+|\lambda|^{4m+2}-|\alpha|^2(|\lambda|^2+|\lambda|^{4m})+(1-|\lambda|^{4m})
(\alpha\bar\lambda-\bar\alpha\lambda)|}\\
&\le\frac{2(1-|\alpha|^2)|\lambda|^{2m+1}}
{1+|\lambda|^{4m+2}-|\alpha|^2(|\lambda|^2+|\lambda|^{4m})}\le
\frac{2|\lambda|^{2m+1}}{1+|\lambda|^{4m+2}},
\end{aligned}
$$ (use that $1+|\lambda|^{4m+2}>|\lambda|^2+|\lambda|^{4m}$).

\begin{proof}[Proof of Proposition \ref{origin}]
Observe that there is a constant $c>0$ with:

\noindent $\bullet$ $c_{A_{m,n}}(p_{m,n}(\lambda),p_{m,n}(\mu))
\ge\max\{\rho(\lambda^m,\mu^m),\rho(\lambda^{m+1},\mu^{m+1})\}
\stackrel{near\ 0}{\ge}c|X_{\lambda,\mu}|;$

\noindent $\bullet$ $\gamma_{A_{m,n}}(0;X_{\lambda,\mu})\ge
c|X_{\lambda,\mu}|;$

\noindent $\bullet$ $|X_{\lambda,\mu}|\ge
c|\lambda^k-\mu^k|\max\{|\lambda|^{k-n},|\mu|^{k-n}|\}$ for any
$k>n$.

Let now $h_{\lambda,\mu}$ be an extremal function for
$c_{A_{m,n}}(p_{m,n}(\lambda),p_{m,n}(\mu)).$ Then $$
h_{\lambda,\mu}(\zeta)=\sum_{j=1}^{[n/m]}a_{j,\lambda,\mu}
\zeta^{jm}+a_{n,\lambda,\mu}\zeta^n+\sum_{j>n,
j\in\S_{m,n}}a_{j,\lambda,\mu}\zeta^j. $$ Since
$|a_{j,\lambda,\mu}|\le 1,$ it follows that $$
|h_{\lambda,\mu}h(\lambda)-h_{\lambda,\mu}(\mu)|\le
H(\lambda,\mu):=$$
$$|a_{m,\lambda,\mu}(\lambda^m-\mu^m)+a_{n,\lambda,\mu}
(\lambda^n-\mu^n)|+\sum_{j=2}^{[n/m]}
|\lambda^{jm}-\mu^{jm}|+\sum_{j=n+1}^\infty|\lambda^j-\mu^j|.$$
Thus, $$1\le\liminf_{\lambda,\mu\to
0,\lambda\neq\mu}\frac{H(\lambda,\mu)}{|h_{\lambda,\mu}(\lambda)-
h_{\lambda,\mu}(\mu)|}=\liminf_{\lambda,\mu\to
0,\lambda\neq\mu}\frac{H(\lambda,\mu)}
{m_{A_{m,n}}(p_{m,n}(\lambda),p_{m,n}(\mu))}$$ $$\le
\liminf_{\lambda,\mu\to
0,\lambda\neq\mu}\frac{|a_{m,\lambda,\mu}(\lambda^m-\mu^m)+a_{n,\lambda,\mu}(\lambda^n-\mu^n)|}
{c_{A_{m,n}}(p_{m,n}(\lambda),p_{m,n}(\mu))}$$ $$+
\liminf_{\lambda,\mu\to 0,\lambda\neq\mu}\frac{\sum_{j=2}^{[n/m]}
|\lambda^{jm}-\mu^{jm}|+\sum_{j=n+1}^\infty|\lambda^j-\mu^j|}{c|X_{\lambda,\mu}|}$$
$$=\liminf_{\lambda,\mu\to
0,\lambda\neq\mu}\frac{|a_m(\lambda^m-\mu^m)+a_n(\lambda^n-\mu^n)|}
{m_{A_{m,n}}(p_{m,n}(\lambda),p_{m,n}(\mu))}
\le\liminf_{\lambda,\mu\to
0,\lambda\neq\mu}\frac{\gamma_{A_{m,n}}(0;X_{\lambda,\mu})}
{c_{A_{m,n}}(p_{m,n}(\lambda),p_{m,n}(\mu))}$$ (since $$
\gamma_{A_{m,n}}(0;X)=\max\{|X_1\frac{h^{(n)}(0)}{n!}+X_2\frac{h^{(m)}(0)}{m!}|:
h\in\CO_{m,n}(\Bbb D),h(0)=0\}).$$

The opposite inequality $$ \limsup_{\lambda,\mu\to
0,\lambda\neq\mu}\frac{\gamma_{A_{m,n}}(0;X_{\lambda,\mu})}
{c_{A_{m,n}}(p_{m,n}(\lambda),p_{m,n}(\mu))}\leq 1 $$ can be
proven in a similar way and we omit the details.
\end{proof}

\begin{proof}[Proof of Corollary \ref{ineq}]
Observe that for any neighborhood $U$ of $0$ we may find points
$\lambda,\mu\in U$ such that $\lambda^m-\mu^m=\lambda^n-\mu^n\neq
0.$ Then, by Proposition \ref{origin}, it is enough to show that
$$\gamma_{A_{m,n}}(0;X_0)>1,\mbox{ where }X_0:=(1,1).$$ Since
$\gamma_{A_{m,n}}(0;X_0)$
$$=\max\{|\frac{h^{(n-m)}(0)}{(n-m)!}+h(0)|:
h\in\CO(\DD,\overline\DD),h^{(j)}(0)=0,j+m\in S_{m,n}\}$$ and
$\ds\max_{s\in S_{m,n}}s=nm-m-n,$ then
$$\gamma_{A_{m,n}}(0;X_0)\ge\max\{|a+b|:(a,b)\in T_{n-m}\},$$
where $T_{n-m}$ is the set of all pairs $(a,b)\in\Bbb C^2$ for
which there is a function $h\in\CO(\DD,\overline\DD)$ of the form
$h(z)=a+bz^{n-m}+o(z^{nm-2m-n}).$

Let $k\in\Bbb N$ be such that $k(n-m)\ge nm-2m-n.$ We shall show
that there is a function $f\in\CO(\DD,\overline\DD)$ of the form
$f(z)=a+bz+o(z^k)$ such that $a,b>0$ and $a+b>1,$ which will imply
that $\gamma_{A_{m,n}}(0;X_0)>1.$

Note that by Shur's theorem (cf. \cite{Gar}) such a function $f$
exists if and only if
$$(1-|a|^2)X_1^2+(1-|a|^2-|b|^2)\sum_{j=2}^nX_j^2\ge2|ab|\sum_{j=2}^n
X_{j-1}X_j,\quad X\in\Bbb R^n.\leqno({3)}$$ Since
$\cos\frac{\pi}{n+1}$ is the maximal eigenvalue of the quadratic
form $\sum_{j=2}^n X_{j-1}X_j,$ it follows that
$$\cos\frac{\pi}{n+1}\sum_{j=1}^nX_j^2\ge\sum_{j=2}^n
X_{j-1}X_j,\quad X\in\Bbb R^n.$$ Then all pairs $(a,b)\in\Bbb C^2$
for which $2\cos\frac{\pi}{n+1}|ab|\le1-|a|^2-|b|^2$ satisfy (3);
in particular, we may choose $a,b>0$ such that $2ab>1-a^2-b^2,$
i.e., $a+b>1.$
\end{proof}

Now we turn to the discussion of the Carath\'eodory-Reiffen
pseudometric on the $(3,4)$-parabola.

\

\noindent{\it Proof of Proposition \ref{34}.} Recall that $$
\gamma_{A_{3,4}}(0;X)=
\max\{|X_1h'(0)+X_2h(0)|:h\in\CO(\DD,\overline\DD), h''(0)=0\}. $$
So, we have to describe the pairs $(a_0,a_1)\in\CC^2$ for which
there is a function $h\in\CO(\DD,\overline\DD)$ of the form
$h(\zeta)=a_0+a_1\zeta+\mbox{o}(\zeta^2)$. Let $I_3$ be the
$3\times3$ unit matrix and $$ M=\left[\begin{matrix} a_0&a_1&0\\
0&a_0&a_1\\ 0&0&a_0\end{matrix}\right]. $$ It follows by Schur's
theorem (cf. \cite{Gar}) that such an $h$ exists if only if
$I_3-M^\ast M$ is a semipositive matrix. It is easy to check that
the last conditions just means that the pair $(|a_0|^2,|a_1|^2)$
belongs to the set $$ C:=\{(a,b)\in\RR_+^2:a+\sqrt b\le
1,ab(1-a)\le((1-a)^2-b)(1-a-b)\}. $$ The second inequality can be
written as $$ b\le(1-a)(1-\sqrt a)\quad \text{ or }\quad
b\ge(1-a)(1+\sqrt a). $$ Hence
$C=\{(a,b)\in\RR_+^2:b\le(1-a)(1-\sqrt a),\ a\le 1\}$. Thus, $$
\begin{aligned}
\gamma_{A_{3,4}}(0;X)&=\max\{|X_1|\sqrt b+|X_2|\sqrt a:(a,b)\in
C\}\\ &=\max\{t\in[0;1]:|X_1|(1-t)\sqrt{1+t}+|X_2|t\}.
\end{aligned}
$$ Straightforward calculations show that the last maximum is
equal to $$
\begin{cases}|X_1|&\text{ if } |X_1|\geq 2|X_2|\\
               |X_2|&\text{ if } |X_2|\geq \sqrt{2}|X_1|\\
               |X_1|\frac{c^3-18c+(c^2+24)^{3/2}}{108}
               &\text{ if
               }1<c:=2\tfrac{|X_2|}{|X_1|}<2\sqrt{2}
\end{cases}.\eqno{\qed}
$$

Now, we shall go to prove Proposition \ref{2k+1}.

\begin{proof}[Proof of Proposition \ref{2k+1}]
Recall that $$
\begin{aligned}
&m_{A_{2,2k+1}}(p_{2,2k+1}(\lambda),p_{2,2k+1}(\mu))\\
&=\max\{m_\DD(f(\lambda),f(\mu)):f\in\CO(\DD,\DD),f^{(2j-1)}(0)=0,j=1,\dots,k\}\\
&{\begin{aligned}=\max\{m_{\DD}(\lambda^2h(\lambda),&
\mu^2h(\mu)):\\ &h\in\CO(\DD,\overline\DD ),
h^{(2j-1)}(0)=0,j=1,\dots,k-1\}.
\end{aligned}}
\end{aligned}
$$ It follows that $$
m_{A_{2,2k+1}}(p_{2,2k+1}(\lambda),p_{2,2k+1}(\mu)) $$ $$
=\sup\{m_{\DD}(\lambda^2z,\mu^2w):m_\DD(z,w) \le
m_{A_{2,2k-1}}(p_{2,2k-1}(\lambda),p_{2,2k-1}(\mu)). $$ Then
Proposition \ref{2k+1} will follow by induction on $n\in\ZZ_+$ if
we show that $$
m_\DD(z,w)\le\frac{2|\lambda|^{2k-1}}{1+|\lambda|^{4k-2}}\implies
m_\DD(\lambda^2z,\lambda^2w)\le\frac{2|\lambda|^{2k+1}}{1+|\lambda|^{4k+2}}.
$$ Since
$\frac{2|\lambda|^{2k-1}}{1+|\lambda|^{4k-2}}=m_\DD(\lambda^{2k-1},-\lambda
^{2k-1})$, we may assume as in Remark (b) that
$z=\Phi_{\alpha}(\lambda^{2k-1})$ and
$w=\Phi_{\alpha}(-\lambda^{2k-1})$ for some $\alpha\in\DD$. Then
$$
\begin{aligned}
&m_{\DD}(\lambda^2z,\lambda^2w)\\
&=\frac{2(1-|\alpha|^2)|\lambda|^{2k+1}}
{|1+|\lambda|^{4k+2}-|\alpha|^2(|\lambda|^4+|\lambda|^{4k-2})+(1-|\lambda|^4)
(\alpha\bar\lambda^{2k-1}-\bar\alpha\lambda^{2k-1})|}\\
&\le\frac{2(1-|\alpha|^2)|\lambda|^{2k+1}}
{1+|\lambda|^{4k+2}-|\alpha|^2(|\lambda|^4+|\lambda|^{4k-2})}\le
\frac{2|\lambda|^{2k+1}}{1+|\lambda|^{4k+2}},
\end{aligned}
$$ since $1+|\lambda|^{4k+2}>|\lambda|^4+|\lambda|^{4k-2}$.
\end{proof}

\noindent{\bf Remark.} From the result above one may conclude the
following interpolation result. Namely, for given $k\in\NN$,
$\lambda, \eta, \zeta\in\DD$ the following conditions are
equivalent:\smallskip

(i) $m_\DD(\eta,\zeta)\leq m_\DD(\lambda^{2k+1},-\lambda^{2k-1});$

(i) $\exists_{f\in\CO(\DD,\DD)}:f(\lambda^{2k+1})=\eta,
f(-\lambda^{2k+1})=\zeta$;

(iii) $\exists_{f\in\CO(\DD,\DD)}:f(\lambda)=\eta,
f(-\lambda)=\zeta, f^{(j)}(0)=0, j=1,\dots ,2k;$

(iv) $\exists_{f\in\CO(\DD,\DD)}:f(\lambda)=\eta,
f(-\lambda)=\zeta, f^{(2j-1)}(0)=0, j=1,\dots ,k.$\smallskip

Indeed, it is trivial that
$(i)\implies(ii)\implies(iii)\implies(iv),$ and the implication
$(iv)\implies(i)$ follows by the equalities
$$m_{A_{2,2k+1}}(p_{2,2k+1}(\lambda),p_{2,2k+1}(\mu))=
\frac{2\lambda|^{2k+1}}{1+|\lambda|^{4k+2}}=m_\DD(\lambda^{2k+1},-\lambda^{2k-1}).$$

Finally, we discuss the proof for the Kobayashi distance and
metric.

\begin{proof}[Proof of Proposition \ref{Koba}]
The proof of the formula for $\wdtl k_{A_{m,n}}$ follows the one
for the case $(m,n)=(2,3)$ (see \cite{Kne}). For convenience of
the reader we include it.

First, $\wdtl
k_{A_{m,n}}(p_{m,n}(\lambda),p_{m,n}(\mu))\le\rho(\lambda,\mu)$
because $p_{m,n}$ is holomorphic. Second, since $m $ and $n$ are
relatively prime, it is easy to see that
$\CO(\DD,A_{m,n})=\{p_{m,n}\circ\psi:\psi\in\CO(\DD,\DD)\}.$ Then
any $\varphi\in\CO(\DD,A_{m,n})$ with
$\varphi(\wdtl\lambda)=p_{m,n}(\lambda)$ and
$\varphi(\wdtl\mu)=p_{m,n}(\mu)$ corresponds to some
$\psi\in\CO(\DD,\DD)$ with $\psi(\wdtl\lambda)=\lambda$ and
$\psi(\wdtl\mu)=\mu.$ Thus, $\rho(\lambda,\mu)\le\rho(\wdtl
\lambda,\wdtl\mu)$ and hence $\rho(\lambda,\mu)\le\wdtl
k_{A_{m,n}}(p_{m,n}(\lambda),p_{m,n}(\mu)).$ So, $\wdtl
k_{A_{m,n}}(p_{m,n}(\lambda),p_{m,n}(\mu))=\rho(\lambda,\mu);$ in
particular, $\wdtl k_{A_{m,n}}$ is a distance and therefore $\wdtl
k_{A_{m,n}}=k_{A_{m,n}}.$

The formulas for $\kappa_{A_{m,n}}$ can be proven in a similar way
and we omit the details.
\end{proof}

We conclude this paper by mentioning the simplest example of a
reducible variety.
\medskip

\noindent{\bf Remark.} Put $A_{2,2}:=\{(z,w)\in\DD^2: z^2=w^2\}$;
$A_{2,2}$ is reducible.  Obviously, $A_{2,2}$ is biholomorphically
equivalent to the coordinate cross $V:=\{(z,w)\in\DD^2:zw=0\}$.
Therefore, we discuss $V$ instead of $A_{2,2}$.

It is clear that $c_V((z_1,0),(z_2,0))=\tilde
k_V((z_1,0),(z_2,0))=\rho(z_1,z_2),$ $$\tilde
k_V((z,0),(0,w))=\infty\ (zw\neq 0)$$ and
$$k_V((z,0),(0,w))=\tilde k_V((z,0),(0,0))+ \tilde
k_V((0,0),(0,w))= \rho(|z|,-|w|).$$ Moreover,
$\gamma_V((z,0);(1,0))=
\kappa_V((z,0);(1,0))=\ds\frac{1}{1-|z|^2}$ and $$
\kappa_V(0;X)=\begin{cases}|X|&\text{ if }X_1X_2=0\\ \infty&\text{
if otherwise}
\end{cases}.
$$

Recall now that $$
\begin{aligned}
\CO(V,\DD) =&\{f+g-f(0):\\ &f\in\CO(\DD\times\{0\},\DD),
g\in\CO(\{0\}\times\DD,\DD),f(0)=g(0)\}.
\end{aligned}
$$ Then obviously $\gamma_V(0;X)=|X_1|+|X_2|.$

Finally, since $z+w\in\CO(V,\DD),$ it follows that
$$c_V((z,0),(0,w))=c_V((|z|,0),(-|w|,0))\ge\rho(|z|,-|w|).$$ Thus,
$c_V=k_V;$ in particular, $c_V=c_V^i.$

\end{document}